\documentclass{amsart}

\usepackage{amsmath}
\usepackage{amssymb}
\usepackage{amsfonts}
\usepackage{amsbsy}
\usepackage{amscd}
\usepackage{eucal}
\usepackage{mathptmx}

\newcommand{\beq}{\begin{equation}}
\newcommand{\eeq}{\end{equation}}

\newtheorem{theorem}{Theorem}
\newtheorem{proposition}{Proposition}

\newtheorem{note}{Note}

\newcommand{\rmd}{\mathrm{d}}
\newcommand{\rmi}{\mathrm{i}}
\newcommand{\Real}{\mathop{\mathrm{Re}}}

\DeclareMathOperator{\sgn}{sgn}

\newcommand{\N}{\mathbb{N}}

\newcommand{\Z}{\mathbb{Z}}
\newcommand{\C}{\mathbb{C}}

\newcommand{\ud}{\frac{1}{2}}
\newcommand{\uds}{{\textstyle\frac{1}{2}}}

\begin{document}

\title[Integral formula for the Bessel function of the first kind]
{Integral formula for the Bessel function of the first kind}

\author[Enrico De Micheli]{Enrico De Micheli}
\address{Consiglio Nazionale delle Ricerche, Via De Marini, 6 - 16149 Genova, Italy}
\email{enrico.demicheli@cnr.it}



\begin{abstract}
In this paper, we prove a new integral formula for the Bessel function of the first kind 
$J_\mu(z)$. This formula generalizes to any $\mu$, $z\in\C$ the classical representations of Bessel and Poisson.
\end{abstract}

\maketitle
 
\section{Introduction}
\label{se:intro}

In this note, we prove the following integral representation of the Bessel function of the first kind 
$J_\mu(z)$, which holds for unrestricted complex values of the order $\mu$:
\beq
J_{\mu}(z) = \frac{1}{2\pi}\left(\frac{z}{2}\right)^\mu
\int_{-\pi}^\pi e^{\rmi z\cos\theta} \, \gamma^*\left(\mu,\uds\rmi z e^{\rmi\theta}\right)\rmd\theta
\qquad (\mu\in\C; z\in\C),
\label{T.0}
\eeq
where $\gamma^*$ denotes the Tricomi version of the (lower) incomplete gamma function. 

Lower and upper incomplete gamma functions arise from the decomposition of the Euler integral for the
gamma function:
\begin{subequations}
\begin{align}
\gamma(\mu,w) &= \int_0^w e^{-t}\,t^{\mu-1}\,\rmd t \qquad (\Real\mu>0), \label{G.1.a} \\
\Gamma(\mu,w) &= \int_w^\infty e^{-t}\,t^{\mu-1}\,\rmd t \qquad (|\arg w|<\pi). \label{G.1.b}
\end{align}
\label{G.1}
\end{subequations}
Analytical continuation with respect to both parameters can be based either on the integrals in 
\eqref{G.1} or on series expansions of $\gamma(\mu,w)$, e.g. \cite[Eq. 9.2(4)]{Erdelyi}:
\beq
\gamma(\mu,w) = \sum_{n=0}^\infty\frac{(-1)^n}{n!}\frac{w^{\mu+n}}{\mu+n},
\label{G.2}
\eeq
which shows that $\gamma(\mu,w)$ has simple poles at $\mu=0,-1,-2,\ldots$ and, in general, is
a multi-valued (when $\mu$ is not a positive integer) function with a branch point at $w=0$.
Similarly, even $\Gamma(\mu,w)$ is in general multi-valued but it is an entire function of $\mu$.
Inconveniences related to poles and multi-valuedness of $\gamma$
can be circumvented by introducing what Tricomi
called the \emph{fundamental function} of the incomplete gamma theory \cite[\S 2, p. 264]{Tricomi2}:
\beq
\gamma^*(\mu,w) = \frac{w^{-\mu}}{\Gamma(\mu)}\int_0^w e^{-t} t^{\mu-1}\,\rmd t
= \frac{e^{-w}}{\Gamma(1+\mu)}\,M(1,1+\mu;w),
\label{G.3}
\eeq
where $M(a,b;z)$ is the Kummer confluent hypergeometric function \cite[Eq. 13.2.2]{DLMF}.
The function $\gamma^*(\mu,w)$ is entire in both $w$ and $\mu$, real-valued for real $w$ and $\mu$
and $\gamma^*(-n,w)=w^n$ for $n=0,1,2,\ldots$ . B\"ohmer \cite{Bohmer} 
and later Tricomi \cite{Tricomi1,Tricomi2} (see also \cite{Gautschi}) 
assumed $\gamma^*(\mu,w)$ as base point 
for the theory, defining the incomplete gamma functions by the relations:
\begin{align}
\gamma(\mu,w) &= \Gamma(\mu) \, w^\mu \, \gamma^*(\mu,w), \label{G.4} \\
\Gamma(\mu,w) &= \Gamma(\mu) \, \left[1-w^\mu\gamma^*(\mu,w)\right]. \label{G.5}
\end{align}
For a detailed account of the incomplete gamma functions theory, see
\cite{Gautschi}, \cite[Chapter 11]{Temme} and \cite{Erdelyi,DLMF}, where massive compilations of results are stated. For what concerns details on the Bessel function, the reader is referred
to the classical treatise of Watson \cite{Watson}.

The strong link between Bessel functions and incomplete gamma functions has long been known.
Examples of this connection are represented, e.g., by the integral representations of $\gamma(\mu,w)$ in terms of $J_\mu(z)$ for $\Real\mu>0$ \cite[Eq. 9.3(4)]{Erdelyi} and of $\Gamma(\mu, w)$ in terms 
of modified Bessel function of the second kind $K_\mu(z)$ for $\Real\mu<1$ \cite[Eq. 8.6.6]{DLMF}.
In this respect, it is worth mentioning also the expansions of incomplete gamma functions in
series of Bessel functions: \cite[Eqs. (2.8), (2.9)]{Gautschi}, \cite[Eqs. (40)]{Tricomi2},
and of modified Bessel functions: \cite[Eqs. 8.7.4, 8.7.5]{DLMF}, \cite[Eqs. (48), (49)]{Tricomi1}.
Formula \eqref{T.0}, whose validity extends to $(\mu,z)\in\C^2$, 
makes even more explicit the intimate link between these two classes of functions.

Classical formulae of the Bessel function of the first kind $J_\mu(z)$ given by integrals 
on the real line can be found in \cite[Chapter 10.9]{DLMF}, \cite[Sect. 7.3]{Erdelyi} and 
\cite[Sect. 8.41]{Gradshteyn}. Many of them hold only for restricted values of $\mu$:
for instance, the classical Bessel integral \cite[p. 70]{Koronev},
\beq
J_\mu(z) = \frac{(-\rmi)^\mu}{\pi}\int_0^\pi e^{\rmi z\cos\theta}\,\cos \mu\theta\,\rmd\theta
\qquad (\mu\in\Z),
\label{T.1}
\eeq
and Poisson's integral \cite[Eq. (14.5), p. 38]{Koronev},
\beq
J_\mu(z) = \frac{1}{\Gamma(\mu+\ud)\sqrt{\pi}}\left(\frac{z}{2}\right)^\mu
\int_0^\pi e^{\rmi z\cos\theta}\,(\sin\theta)^{2\mu}\,\rmd\theta
\qquad (\Real\mu>-\uds),
\label{T.2}
\eeq
hold only for $\mu\in Z$ and $\Real\mu>-\ud$, respectively.
Additional examples are the Gubler representations \cite[Eq. 7.3(11)]{Erdelyi},
which holds for $\Real\mu<\ud$, the Mehler-Sonine integral formulae for $|\Real\mu|<\ud$ 
\cite[Eq. 7.12(12)]{Erdelyi} or for $|\Real\mu|<1$ \cite[Eq. 7.12(14)]{Erdelyi}. 
Only the well-known Schl\"afli's and Heine's 
representations hold for unrestricted complex values of the index $\mu$. Schl\"afli's 
representation \cite[Eq. 7.3(9)]{Erdelyi} holds for $\mu\in\C$, with $\Real z>0$ 
(it holds also in case $\Real z=0$ provided that $\Real\mu>0$),
while the two Heine's expressions \cite[Eqs. 7.3(31) \& (32)]{Erdelyi} hold for $\mu\in\C$, but separately in the upper and lower half-plane of the complex $z$-plane, respectively. 
In both cases, Schl\"afli's and Heine's, the representations are made up of the sum of two integrals.

Standard generalizations of Poisson's integral are obtained via its extension 
to complex contour integrals, examples of which are the classical Hankel's representations 
\cite[Eqs. (19.7), (19.8)]{Koronev}.
As Poisson's representation of $J_\mu(z)$ can be seen as the generalization of Bessel's integral
to complex values of $\mu$ with $\Real\mu>-\ud$ (recall that the trivial extension of formula 
\eqref{T.1} to complex values of $\mu$ yields the Anger function $\mathbf{J}_\mu(z)$, instead), 
mere inspection shows that formula 
\eqref{T.0} can be viewed as the generalization of Poisson's integral to represent $J_\mu(z)$ 
at values of $\mu$ belonging to the entire complex plane, the role of $(\sin\theta)^{2\mu}$ in \eqref{T.2} being now taken by the incomplete gamma function $\gamma^*(\mu,\ud\rmi ze^{\rmi\theta})$.

The proof of formula \eqref{T.0} is given in the next section and, essentially, makes use
of elementary properties of the special functions involved. A preliminary result was given in 
\cite{DeMicheli}.

\section{Integral representation for $J_\mu(z)$}
\label{se:J}

We denote by $\N_0=\{0,1,2,\ldots\}$ the set of non-negative integers.

\begin{theorem}
The following integral representation of the Bessel function of the first kind
holds true for any complex order $\mu\in\C$ and for $z\in\C$ (slit along the real negative axis when $\mu\not\in\Z$):
\beq
J_{\mu}(z) = \frac{1}{2\pi}\left(\frac{z}{2}\right)^\mu
\int_{-\pi}^\pi e^{\rmi z\cos\theta} \, \gamma^*\left(\mu,\uds\rmi z e^{\rmi\theta}\right)\rmd\theta
\qquad (\mu\in\C; z\in\C).
\label{T.00}
\eeq
\end{theorem}

\begin{proof}
Our starting point is the Gegenbauer generalization of Poisson's integral representation of the 
Bessel functions of the first kind \cite[\S 3.32, Eq. (1), p. 50]{Watson}, which holds for
$\ell\in\N_0$ and $\Real\nu>-\uds$:
\beq
J_{\nu+\ell}(z) = \frac{(-\rmi)^\ell\Gamma(2\nu)\,\ell!}{\Gamma(\nu+\uds)\Gamma(\uds)\Gamma(2\nu+\ell)}
\left(\frac{z}{2}\right)^\nu\int_0^\pi e^{\rmi z\cos u}(\sin u)^{2\nu} 
C_\ell^\nu(\cos u)\,\rmd u,
\label{d.J.4}
\eeq
where $C_\ell^\nu(t)$ is the Gegenbauer polynomial of order $\nu$ and degree $\ell$.
Now, we recall from \cite[Proposition 1]{DeMicheli1} the following representation of the 
Gegenbauer polynomials $C_\ell^\nu(\cos u)$, which holds for $\Real\nu>0$ and $\ell\in\N_0$:
\beq
C_\ell^{\nu}(\cos u) = \frac{\Gamma(\ell+2\nu)}{2^\nu\,\ell!\,[\Gamma(\nu)]^2}
\frac{e^{-\rmi\pi\nu}}{(\sin u)^{(2\nu-1)}}
\int_u^{2\pi-u} \!\! e^{\rmi(\ell+\nu)t} \ (\cos u-\cos t)^{\nu-1}\,\rmd t \quad (u\in[0,\pi]).
\label{d.22}
\eeq
Plugging \eqref{d.22} into \eqref{d.J.4} and using the Legendre duplication formula
for the gamma function, we obtain:
\beq
J_{\nu+\ell}(z) = \frac{(-\rmi)^\ell\,e^{-\rmi\pi\nu}\, z^\nu}{2\pi\,\Gamma(\nu)}
\int_0^\pi\!\rmd u\,\sin u\,e^{\rmi z\cos u}
\int_u^{2\pi-u}\!\!e^{\rmi(\ell+\nu)\theta}\,(\cos u-\cos\theta)^{\nu-1}\,\rmd\theta.
\label{d.J.6}
\eeq
Interchanging the order of integration, \eqref{d.J.6} can be written as follows:
\beq
\begin{split}
J_{\nu+\ell}(z) &= \frac{(-\rmi)^\ell\,e^{-\rmi\pi\nu}z^\nu}{2\pi\,\Gamma(\nu)}
\left[\int_{0}^\pi\!\rmd\theta\,e^{\rmi(\ell+\nu)\theta}\int_{0}^\theta \!
e^{\rmi z\cos u} (\cos u-\cos\theta)^{\nu-1}\,\sin u\,\rmd u \right. \\
&\quad\left. +\int_\pi^{2\pi}\!\rmd\theta\,e^{\rmi(\ell+\nu)\theta}
\int_{0}^{2\pi-\theta}\!e^{\rmi z\cos u}\,(\cos u-\cos\theta)^{\nu-1}\,\sin u\,\rmd u \right].
\end{split}
\label{d.J.7}
\eeq
Next, changing the integration variables: $\theta-2\pi\to\theta$ and $u\to-u$,
the second integral on the r.h.s. of \eqref{d.J.7} becomes:
\beq
e^{\rmi\nu 2\pi}\int_{-\pi}^0 \rmd\theta \ e^{\,\rmi(\ell+\nu)\theta}
\int_0^\theta e^{\rmi z\cos u}(\cos u-\cos\theta)^{\nu-1} \sin u \,\rmd u,
\label{d.J.7.1} \nonumber
\eeq
which, inserted in \eqref{d.J.7}, yields:
\beq
\begin{split}
J_{\nu+\ell}(z) &= \frac{(-\rmi)^\ell z^\nu}{2\pi\,\Gamma(\nu)}
\left[e^{\rmi\nu\pi}\int_{-\pi}^0\!\rmd\theta\,e^{\,\rmi(\ell+\nu)\theta}
\int_{0}^\theta\!e^{\rmi z\cos u} (\cos u-\cos\theta)^{\nu-1}\sin u\,\rmd u \right. \\
&\quad\left. +e^{-\rmi\nu\pi}\int_0^{\pi}\!\rmd\theta\,e^{\rmi(\ell+\nu)\theta}\int_{0}^{\theta} \! 
e^{\rmi z\cos u} (\cos u-\cos\theta)^{\nu-1}\sin u\,\rmd u \right].
\end{split}
\label{d.J.8}
\eeq
Formula \eqref{d.J.8} allows us to write $J_{\nu+\ell}(z)$ as the Fourier coefficient
of a suitable $2\pi$-periodic function, i.e.:
\beq 
\rmi^\ell J_{\nu+\ell}(z) = \widetilde{\CMcal{A}}(\nu,z) \doteq
\int_{-\pi}^\pi \CMcal{A}(\theta;\nu,z) \ e^{\rmi\ell\theta}\,\rmd\theta 
\qquad (\ell\in\N_0, \Real\nu>0),
\label{d.J.9}
\eeq
where $\CMcal{A}(\theta;\nu,z)$ denotes the Abel-type integral
\beq
\begin{split}
\CMcal{A}(\theta;\nu,z) &= 
\frac{z^\nu}{2\pi\Gamma(\nu)}e^{\rmi\nu[\theta-\pi\sgn(\theta)]} 
\int_{\cos\theta}^1 e^{\rmi zt}(t-\cos\theta)^{\nu-1}\,\rmd t \qquad (\Real\nu>0),
\end{split}
\label{d.J.10}
\eeq
and $\sgn(\cdot)$ denotes the \emph{sign function}.
Changing in \eqref{d.J.10} the integration variable $-\rmi z(t-\cos\theta)\to t$
and recalling \eqref{G.1.a}, $\CMcal{A}(\theta;\nu,z)$ can be rewritten as follows:
\beq
\CMcal{A}(\theta;\nu,z) = \frac{\rmi^\nu}{2\pi} 
e^{\rmi\nu[\theta-\pi\sgn(\theta)]}e^{\rmi z\cos\theta} 
\ P(\nu,-\rmi z\,(1-\cos\theta))  \qquad (\Real\nu>0),
\label{d.J.10.bis}
\eeq
where $P(\nu,w)\doteq\gamma(\nu,w)/\Gamma(\nu)$ denotes the normalised (lower) incomplete gamma function. Finally, with some algebra, we obtain from \eqref{d.J.9} and \eqref{d.J.10.bis}
the representation
\beq
J_{\ell+\nu}(z) = \frac{\rmi^{(\ell+\nu)}}{2\pi}
\int_{-\pi}^\pi e^{\rmi z\cos\theta} \, e^{-\rmi(\ell+\nu)\pi\sgn(\theta)}
\ P(\nu,-\rmi z\,(1-\cos\theta)) \ e^{\rmi(\ell+\nu)\theta}\,\rmd\theta,
\label{f.3}
\eeq
which holds for $\ell\in\N_0$ and $\Real\nu>0$.
Now, we want to show that the integral in \eqref{f.3} does not change if we substitute the factor
$P(\nu,-\rmi z\,(1-\cos\theta))$ with
$P(\ell+\nu,-\rmi z\,(1-\cos\theta))$. To this end, we recall that $P(\nu,w)$ satisfies
the following recurrence relation \cite[Eq. 8.8.11]{DLMF}:
\beq
P(\ell+\nu,w)=P(\nu,w)-w^\nu\,e^{-w}\sum_{k=0}^{\ell-1}\frac{w^k}{\Gamma(\nu+k+1)}
\qquad (\ell\geqslant 0),
\label{f.3.1}
\eeq
the sum being obviously understood to be null when $\ell=0$.
We have:
\beq
\begin{split}
&\int_{-\pi}^\pi e^{\rmi z\cos\theta} \, e^{-\rmi(\ell+\nu)\pi\sgn(\theta)}
\ P(\ell+\nu,-\rmi z\,(1-\cos\theta)) \ e^{\rmi(\ell+\nu)\theta}\,\rmd\theta \\
&=\int_{-\pi}^\pi e^{\rmi z\cos\theta} \, e^{-\rmi(\ell+\nu)\pi\sgn(\theta)}
\ P(\nu,-\rmi z\,(1-\cos\theta)) \ e^{\rmi(\ell+\nu)\theta}\,\rmd\theta \\
&\quad-(-\rmi z)^\nu e^{\rmi z} \sum_{k=0}^{\ell-1}\frac{(-\rmi z)^k}{\Gamma(\nu+k+1)}
\int_{-\pi}^\pi e^{-\rmi(\ell+\nu)\pi\sgn(\theta)}(1-\cos\theta)^{\nu+k}e^{\rmi(\ell+\nu)\theta}
\,\rmd\theta.
\end{split}
\label{f.3.2}
\eeq
By using formula \cite[Formula 3.63.9]{Gradshteyn}, the rightmost integrals on the r.h.s. 
of \eqref{f.3.2} can be readily computed:
\beq
\begin{split}
&\int_{-\pi}^\pi e^{-\rmi(\ell+\nu)\pi\sgn(\theta)}(1-\cos\theta)^{\nu+k}e^{\rmi(\ell+\nu)\theta}
\,\rmd\theta
= 2^{\nu+k+2} \int_{0}^{\pi/2} (\cos\theta)^{2(\nu+k)}\,\cos[2(\ell+\nu)\theta]\,\rmd\theta \\
&=\frac{2^{1-\nu-k}\pi}{(2\nu+2k+1) \ B(2\nu+k+\ell+1,k-\ell+1)},
\label{f.3.3}
\end{split}
\eeq
where $B(\cdot,\cdot)$ denotes the Euler beta function. For $\Real\nu>-\ud$ and $\ell\geqslant 0$, 
the integrals in \eqref{f.3.3} are therefore null for $0\leqslant k\leqslant \ell-1$.  
Hence, in view of this latter result, we now set $\mu=\ell+\nu$, and from \eqref{f.3} we have:
\beq
J_{\mu}(z) = \frac{\rmi^{\mu}}{2\pi}
\int_{-\pi}^\pi e^{\rmi z\cos\theta} \, e^{-\rmi\mu\pi\sgn(\theta)}
\ P(\mu,-\rmi z\,(1-\cos\theta)) \ e^{\rmi\mu\theta}\,\rmd\theta
\qquad (\Real\mu>0).
\label{f.3.4}
\eeq
Now, we change the integration variable $\theta-\pi\sgn(\theta) \to \theta$, write
$(1+\cos\theta) = \ud e^{-\rmi\theta}(1+e^{\rmi\theta})^2$ and, in view of the 
$2\pi$-periodicity of the integrand, formula \eqref{f.3.4} becomes:
\beq
J_{\mu}(z) = \frac{\rmi^{\mu}}{2\pi}
\int_{-\pi}^\pi e^{-\rmi z\cos\theta} \, P(\mu,-\uds\rmi z e^{-\rmi\theta}(1+e^{\rmi\theta})^2) 
\, e^{\rmi\mu\theta}\,\rmd\theta \qquad (\Real\mu>0).
\label{f.4}
\eeq
Note that $P(\mu,w)$ is an entire function of $\mu$ but it is a multi-valued function of $w$
with a branch point at $w=0$.

\vspace{0.2cm}

Consider for $\lambda\in[0,\infty)$ the integral
\beq
S_{z,\mu}(\lambda)\doteq
\int_{-\pi}^\pi s_{z,\mu}(\theta,\lambda) \,\rmd\theta \qquad (\Real\mu>0,z\in\C; \lambda\in[0,+\infty)),
\label{d.1}
\eeq
where
\beq
s_{z,\mu}(\theta,\lambda) \doteq e^{-\rmi z\cos\theta} \, P(\mu,g_\theta(\lambda)) 
\, e^{\rmi\mu\theta} \qquad (\theta\in(-\pi,\pi), \lambda\in[0,+\infty);\Real\mu>0),
\label{d.1.1}
\eeq
with $g_\theta(\lambda)=-\uds\rmi z e^{-\rmi\theta}(1+e^{\rmi\theta})^\lambda$.
Evidently, $2\pi\,J_\mu(z)=\rmi^\mu S_{z,\mu}(2)$.
The function $S_{z,\mu}(\lambda)$ is actually independent of $\lambda$ for $\lambda\geqslant0$
and $\Real\mu>0$. To see this, we analyze the derivative with respect to $\lambda$, 
$\partial_\lambda S_{\mu,z}(\lambda)$.
The function $s_{z,\mu}(\theta,\lambda)$ is continuous in 
$(\theta,\lambda) \in (-\pi,\pi) \times [0,\infty)$ with
\beq
\int_{-\pi}^\pi \left|s_{z,\mu}(\theta,\lambda)\right|\,\rmd\theta
\leqslant 
c_1(\mu,z) \, \exp\left[2^{\lambda-1}|z|+(\lambda-1)|\mu|\right] < \infty \qquad (\Real\mu>0),
\label{d.2.1}
\eeq
for finite values of $\lambda$, $c_1(\mu,z)$ being a finite constant indipendent of $\lambda$.
Similarly, recalling that $\partial_w P(\mu,w) = w^{\mu-1}\,e^{-w}/\Gamma(\mu)$ \cite[Eq. 8.8.13]{DLMF}, we have
\beq
\partial_\lambda s_{z,\mu}(\theta,\lambda)
= \frac{e^{-\rmi z\cos\theta}}{\Gamma(\mu)}\,[g_\theta(\lambda)]^\mu\,e^{-g_\theta(\lambda)}\,\ln(1+e^{\rmi\theta})\,e^{\rmi\mu\theta},
\label{d.2.2}
\eeq
which is continuous in $(\theta,\lambda) \in (-\pi,\pi) \times [0,\infty)$ with:
\beq
\int_{-\pi}^\pi \left|\partial_\lambda s_{z,\mu}(\theta,\lambda)\right|\,\rmd\theta
\leqslant 
c_2(\mu,z) \, \exp[2^{\lambda-1}|z|+(\lambda-1)|\mu|]\int_{-\pi}^\pi\left|\ln(1+e^{\rmi\theta})\right|\,\rmd\theta < \infty
\label{d.2.3}
\eeq
for finite values of $\lambda$, $c_2(\mu,z)$ being a finite constant indipendent of $\lambda$. 
Hence, differentiation under the integral sign is legitimate and we have:
$\partial_\lambda S_{z,\mu}(\lambda)= \int_{-\pi}^\pi\partial_\lambda s_{z,\mu}(\theta,\lambda)\,\rmd\theta$ for $\lambda\geqslant 0$.

Let us introduce the complex variable $w=1+e^{\rmi\theta}$ and, for $\Real\mu>0$ and
$\lambda\geqslant 0$, the function
\beq
f_{z,\mu,\lambda}(w) \doteq \frac{1}{\Gamma(\mu)}\left(-\frac{\rmi z}{2}\right)^\mu
\exp\left(-\frac{\rmi z}{2}\frac{w^2-w^\lambda-2w+2}{w-1}\right)
w^{\lambda\mu}\,\ln w \qquad (\Real\mu>0,\lambda\geqslant 0).
\label{d.6}
\eeq
We have:
\beq
\partial_\lambda S_{\mu,z}(\lambda)=
\oint_{C(1,1)} f_{z,\mu,\lambda}(w) \, \rmd w,
\label{d.7}
\eeq
where $C(1,1)$ denotes the circle with center $w=1$ and radius $R=1$.
The points of $f_{z,\mu,\lambda}(w)$ which require attention are $w=0$ and $w=1$.
We have (recall that $\Real\mu>0$):
\beq
\lim_{w\to 0} f_{z,\mu,\lambda}(w) = \left(-\frac{\rmi z}{2}\right)^\mu \frac{e^{\rmi z}}{\Gamma(\mu)}
\lim_{w\to 0} w^{\lambda\mu}\ln w = 0 \quad\mathrm{for}\quad\lambda>0.
\label{d.8}
\eeq
Moreover, $\lim_{\lambda\to 0^+}\lim_{w\to 0} f_{z,\mu,\lambda}(w) = 0$,
the two limits being evidently not interchageable. 
For what concerns the point $w=1$, we have:
\beq
\lim_{w\to 1} f_{z,\mu,\lambda}(w) = \frac{1}{\Gamma(\mu)}
\left(-\frac{\rmi z}{2}\right)^\mu e^{\rmi \lambda z/2} \lim_{w\to 1} \ln w = 0
\quad\mathrm{for}\quad\lambda\in[0,\infty).
\label{d.10}
\eeq
Hence, for any $\lambda>0$, $f_{z,\mu,\lambda}(w)$ is a single-valued holomorphic function 
with no singularities on and within the circle $C(1,1)$. Then, by Cauchy theorem the integral in
\eqref{d.7} vanishes and $\partial_\lambda S_{\mu,z}(\lambda) = 0$ for $\lambda>0$.
Moreover, from \eqref{d.2.2} we have for $\lambda=0$:
\beq
\left.\partial_\lambda S_{\mu,z}(\lambda)\right|_{\lambda=0}
=\frac{1}{\Gamma(\mu)}\left(-\frac{\rmi z}{2}\right)
\int_{-\pi}^\pi \exp\left[-\rmi z e^{\rmi\theta}/2\right]\ln(1+e^{\rmi\theta})\,\rmd\theta=0
\label{d.10.1}
\eeq
For the continuity of $\partial_\lambda S_{\mu,z}(\lambda)$ for $\lambda\geqslant 0$ which follows
from \eqref{d.2.2} and \eqref{d.2.3}, $\partial_\lambda S_{\mu,z}(\lambda) = 0$ 
for $\lambda\geqslant 0$ and $\Real\mu>0$.
Consequently, from \eqref{f.4} we can write:
\beq
J_{\mu}(z) = \frac{\rmi^{\mu}}{2\pi}
\int_{-\pi}^\pi e^{-\rmi z\cos\theta} \, P(\mu,-\uds\rmi z e^{-\rmi\theta})\,e^{\rmi\mu\theta}\,\rmd\theta 
\qquad (\Real\mu>0).
\label{f.5}
\eeq
Since $P(\mu,w)=w^\mu\,\gamma^*(\mu,w)$ (see \eqref{G.4}), 
\eqref{f.5} yields:
\beq
J_{\mu}(z) = \frac{1}{2\pi}\left(\frac{z}{2}\right)^\mu
\int_{-\pi}^\pi e^{-\rmi z\cos\theta} \, \gamma^*(\mu,-\uds\rmi z e^{-\rmi\theta})\,\rmd\theta
\qquad (\Real\mu>0),
\label{f.6}
\eeq
which, simply changing the integration variable $\theta\to\pi-\theta$ and by the 
$2\pi$-periodicity of the integrand, finally gives formula \eqref{T.0}.
The restriction $\Real\mu>0$ is immaterial since both parts of the equality \eqref{f.6}
are analytic functions of $\mu$ (recall that $\gamma^*(\mu,w)$ is entire in $\mu$ and $w$)
and, according to the principle of analytic continuation, the result is valid for any $\mu\in\C$
and $z\in\C$ (this latter domain being slit along the real negative axis when $\mu$ is not integer).
\end{proof}

It is obvious that, by a trivial change of variable, formula \eqref{T.0} can be rewritten in terms of
the generating function for the Bessel functions of the first kind of integral order:
\beq
e^{\rmi z\sin\theta} = \sum_{n=-\infty}^\infty J_n(z) \, e^{\rmi nz}.
\label{Generating}
\eeq
We have:
\beq
J_\mu(z) = \frac{1}{2\pi}\left(\frac{z}{2}\right)^\mu \int_{-\pi}^\pi
e^{\rmi z\sin\theta}\,\gamma^*\left(\mu,\uds z e^{\rmi\theta}\right)\,\rmd\theta
\qquad (\mu\in\C,z\in\C).
\label{Generating.1}
\eeq
Formula \eqref{T.0} can be extended in following way. First, recall the following recurrence 
relation for $\gamma^*(\mu,w)$, which follows from \eqref{f.3.1} and \eqref{G.4}:
\beq
\gamma^*(\mu+n,w)
=w^{-n}\left(\gamma^*(\mu,w)-e^{-w}\sum_{k=0}^{n-1}\frac{w^k}{\Gamma(\mu+k+1)}\right)
\qquad (n\geqslant 0).
\label{f.7}
\eeq
Then we consider $J_{\mu+n}(z)$ with $n\geqslant 0$. From \eqref{T.0} and \eqref{f.7} we have:
\beq
\begin{split}
J_{\mu+n}(z) &= \frac{1}{2\pi}\left(\frac{z}{2}\right)^\mu
\int_{-\pi}^\pi e^{\rmi z\sin\theta} \, \gamma^*\left(\mu,\uds z e^{\rmi\theta}\right) 
e^{-\rmi n\theta}\,\rmd\theta \\
& \quad-\sum_{k=0}^{n-1}\frac{(z/2)^{\mu+k}}{\Gamma(\mu+k+1)}
\int_{-\pi}^\pi \exp(-\uds ze^{-\rmi\theta}) e^{\rmi(k-n)\theta}\,\rmd\theta
\quad (n\geqslant 0).
\end{split}
\label{f.8}
\eeq
Now, it is easy to show that the rightmost integrals on the r.h.s. of \eqref{f.8} are
null for $k \leqslant n-1$. Therefore, we have:
\beq
J_{\mu+n}(z) = \frac{1}{2\pi}\left(\frac{z}{2}\right)^\mu
\int_{-\pi}^\pi e^{\rmi z\sin\theta} \, \gamma^*(\mu,\uds z e^{\rmi\theta}) 
e^{-\rmi n\theta}\,\rmd\theta \qquad (n\geqslant 0;\mu\in\C,z\in\C),
\label{f.9}
\eeq
that is, for $n\geqslant0$, $2\pi(z/2)^{-\mu}J_{\mu+n}(z)$ coincides with the 
Fourier coefficient $\widetilde{\kappa}_{-n}(\mu,z)$ (see \eqref{d.J.9}) of the function 
$\kappa(\theta;\mu,z)\doteq e^{\rmi z\sin\theta} \, \gamma^*(\mu,\uds z e^{\rmi\theta})$.
Representation \eqref{Generating.1} can thus be viewed as the particular case 
with $n=0$ of formula \eqref{f.9}.  Analogously, the Fourier coefficients 
$\widetilde{\kappa}_{n}(\mu,z)$ for $n>0$ are easily proved to be
\beq
\widetilde{\kappa}_{n}(\mu,z)
= 2\pi\left(\frac{z}{2}\right)^{-\mu}\sum_{k=n}^\infty\frac{(-1)^k}{k!}\frac{1}{\Gamma(\mu-n+k+1)}
\left(\frac{z}{2}\right)^{2k+\mu-n} \quad (n > 0),
\label{fourier.1}
\eeq
the sum on the r.h.s. of \eqref{fourier.1} being the \emph{shifted} 
(i.e., starting from $k=n$) power series defining $J_{\mu-n}(z)$ (see \cite[Eq. 10.2.2]{DLMF}).

\vspace{0.5cm}

The structure of formula \eqref{T.0} (and \eqref{Generating.1}) yields readily analogous
companion results for functions related to the Bessel function of the first kind. 
To exemplify, a couple of instances are given below.

\begin{proposition}
\label{pro:1}
The modified Bessel function of the first kind $I_\mu(z)$ has the following integral representation
for unrestricted complex values of $\mu$ and $z$: 
\beq
I_\mu(z) = \frac{1}{2\pi}\left(\frac{z}{2}\right)^{\mu}
\int_{-\pi}^\pi e^{z\cos\theta}\,\gamma^*(\mu,\uds z e^{\rmi\theta})\,\rmd\theta
\qquad (\mu\in\C; -\pi\leqslant\arg z<\pi).
\label{T.4}
\eeq
\end{proposition}
\begin{proof}
Formula \eqref{T.4} follows from the connection relation \cite[Eq. 10.27.6]{DLMF}:
\beq
I_\mu(z) = e^{\mp\rmi\mu\pi/2}J_\mu(z e^{\pm\rmi\pi/2})
\qquad \left(-\pi\leqslant\pm\arg z\leqslant{\frac{\pi}{2}}\right).
\label{T.3}
\eeq
\end{proof}
As an additional example, from formula \eqref{Generating.1}, an analogous new integral representation for the derivatives of $J_\mu(z)$
follows.

\begin{proposition}
\label{pro:2}
For any $\mu,z\in\C$, the $k$th derivative with respect to $z$ of the Bessel function of 
the first kind, $\partial^k_z J_\mu(z)$, ($k\in\N_0$), can be represented by the following integral:
\beq
\partial_z^k J_\mu(z)
%\frac{\partial^k J_\mu(z)}{\partial z^k} 
= \frac{1}{2\pi}\left(\frac{z}{2}\right)^{\mu-k}
\int_{-\pi}^\pi e^{\rmi z\sin\theta}\,\gamma^*(\mu-k,\uds z e^{\rmi\theta})(\rmi\sin\theta)^k
\,e^{-\rmi k\theta}\,\rmd\theta \quad (k\in\N_0).
\label{T.10}
\eeq
\end{proposition}

\begin{proof}
Representation \eqref{T.10} follows immediately by plugging \eqref{Generating.1} into the 
following formula for the derivative of the Bessel function of the first kind \cite[Eq. 10.6.7]{DLMF}:
\beq
%\frac{\partial^k J_\mu(z)}{\partial z^k} 
\partial_z^k J_\mu(z) = 2^{-k}\sum_{m=0}^k (-1)^m\binom{k}{m}\,J_{\mu-k+2m}(z),
\label{T.11}
\eeq
and recalling that $\sum_{m=0}^k (-1)^m\binom{k}{m}\exp(-2\rmi m\theta)=(1-\exp(-2\rmi\theta))^k$.
\end{proof}

Formula \eqref{T.10} can be analytically continued to $k\in\C$ with $\Real k >-1$, thus yielding an integral representation for the fractional derivative of the Bessel function of the first kind.

\begin{note}
\rm
In view of relation \eqref{G.3} between the Tricomi incomplete gamma function $\gamma^*$ and the Kummer confluent hypergeometric function $M$, the main results of this work can be formulated
also in terms of the latter. Therefore, the integral representation of the Bessel function of the first
kind \eqref{T.00} reads:
\beq
J_{\mu}(z) = \frac{1}{2\pi}\frac{(z/2)^\mu}{\Gamma(1+\mu)}
\int_{-\pi}^\pi e^{\ud\rmi z \, e^{-\rmi\theta}} \, M\left(1,1+\mu,\uds\rmi z e^{\rmi\theta}\right)\rmd\theta
\qquad (\mu\in\C; z\in\C).
\label{T.12}
\eeq
The modified Bessel function of the first order $I_\mu(z)$ reads (see \eqref{T.4}):
\beq
I_\mu(z) = \frac{1}{2\pi}\frac{(z/2)^\mu}{\Gamma(1+\mu)}
\int_{-\pi}^\pi e^{\ud z \, e^{-\rmi\theta}} \,M\left(1,1+\mu;\uds z e^{\rmi\theta}\right)\,\rmd\theta
\qquad (\mu\in\C; -\pi\leqslant\arg z<\pi),
\label{T.13}
\eeq
while the integral representation of the $z$-derivative of fractional order $k$, with $k\in\C$ and $\Real k >-1$, of the Bessel function of the first kind, $\partial_z^k J_\mu(z)$, can be written as
(see\eqref{T.10}):
\beq
\partial_z^k J_\mu(z)
= \frac{1}{2\pi}\frac{(z/2)^{\mu-k}}{\Gamma(1+\mu-k)}
\int_{-\pi}^\pi e^{-\ud z\, e^{-\rmi\theta}}\,M\left(1,1+\mu-k,\uds z e^{\rmi\theta}\right)(\rmi\sin\theta)^k
\,e^{-\rmi k\theta}\,\rmd\theta,
\label{T.14}
\eeq
where the principal branch of the power is taken when $k$ is a non-integer number.
\end{note}

\bibliographystyle{amsplain}

\end{document}